\documentclass[11pt,twoside]{amsart}

\usepackage[dvipsnames]{xcolor}
\usepackage{hyperref}
\hypersetup{colorlinks=true, linkcolor=Blue, citecolor=Maroon}
\usepackage{amsthm, amsmath, amscd, amssymb,centernot,txfonts}
\usepackage{tikz-cd}
\usepackage{tikz}

\usepackage{lipsum}
\usepackage{multicol}
\usepackage{amsfonts}
\usepackage{adjustbox}
\usepackage{mathrsfs}
\usepackage[all]{xy}
\usepackage[left=2.5cm,top=2.5cm,bottom=3cm,right=2.5cm]{geometry}
\usepackage{dirtytalk}
\usepackage{mathtools}
\usepackage{enumitem}
\usepackage{comment}
\usepackage{fourier} 
\usepackage{capt-of}
\usetikzlibrary{shapes,arrows}

\theoremstyle{plain}
\numberwithin{equation}{section}
\newtheorem{theorem}{Theorem}[section]
\newtheorem{proposition}[theorem]{Proposition}
\newtheorem{lemma}[theorem]{Lemma}
\newtheorem{corollary}[theorem]{Corollary}

\theoremstyle{definition}

\newtheorem{definition}[theorem]{Definition}

\newtheorem{set-up}[theorem]{Set-up}

\newcommand*{\QEDB}{\hfill\ensuremath{\square}}

\tikzstyle{decision} = [diamond, draw, , 
    text width=4.5em, text badly centered, node distance=3cm, inner sep=0pt]
\tikzstyle{block} = [rectangle, draw, , 
    text width=10em, text centered, rounded corners, minimum height=2em]
\tikzstyle{block1} = [rectangle, draw, , 
    text width=5em, text centered, rounded corners, minimum height=2em]    
\tikzstyle{line} = [draw, -latex']
\tikzstyle{cloud} = [draw, ellipse,, node distance=3cm,
    minimum height=2em]
\definecolor{light-gray}{gray}{0.95}    
\begin{document}

\title[A note on stability of syzygy bundles on Enriques and bielliptic surfaces]{A note on stability of syzygy bundles on Enriques and bielliptic surfaces}

\author[J. Mukherjee]{Jayan Mukherjee}
\address{The Institute for Computational and Experimental Research in Mathematics, Providence, USA}
\email{jayan\textunderscore mukherjee@brown.edu}

\author[D. Raychaudhury]{Debaditya Raychaudhury}
\address{The Fields Institute for Research in Mathematical Sciences, Toronto, Canada}
\email{draychau@fields.utoronto.ca}


\maketitle
\begin{abstract} In this note, we prove that the syzygy bundle $M_L$ is cohomologically stable with respect to $L$ for any ample and globally generated line bundle $L$ on an Enriques (resp. bielliptic) surface over an algebraically closed field of characteristic $\neq 2$ (resp. $\neq 2,3$). In particular our result on complex Enriques surfaces improves the result of \cite{TLZ21}, Corollary 3.5 by removing the condition on Clifford index. Together with the results of Camere (\cite{Cam12}) and Caucci--Lahoz (\cite{CL21}), it implies that $M_L$ is stable with respect to $L$ for an ample and globally generated line bundle $L$ on any smooth minimal complex projective surface $X$ of Kodaira dimension zero.
\end{abstract}

\section{Introduction} Let $X$ be a smooth projective variety over an algebraically closed field $\mathbb{k}$ and let $L$ be a globally generated line bundle on $X$. The {\it syzygy bundle} $M_L$ associated to $L$ is defined as the kernel of the evaluation map of the global sections of $L$ i.e., by the following short exact sequence 
\begin{equation}\label{*}
    0\to M_L\to H^0(L)\otimes \mathcal{O}_X\to L\to 0.
\end{equation}

One can similarly define the syzygy bundle $M_V$ associated to a base point free sublinear series $V\subseteq H^0(L)$ as the kernel of the evaluation map $V\otimes \mathcal{O}_X\to L$. The study of stability of syzygy bundles is an active research area and has been investigated by many authors. Ein and Lazarsfeld showed that for a smooth projective curve $C$ of genus $g$ over an algebraically closed field of arbitrary characteristic, and a globally generated line bundle $L$ on $C$, the syzygy bundle $M_L$ is stable as soon as $\textrm{deg}(L)\geq 2g+1$ (\cite{EL92}, Proposition 1.5). The stability of syzygy bundles for projective spaces (both for complete and incomplete linear series) was investigated extensively by Flenner (\cite{Fle84}), Trivedi (\cite{Tri10}), Brenner (\cite{Bre08}), Costa, Marcias Marques, Mir\'o-Roig (\cite{CMM10}, \cite{MM11}) and Coand\u{a} (\cite{Coa11}) among others.  

The main motivation of this note comes from the works of Camere, Ein--Lazarsfeld--Mustopa and Torres-L\'opez--Zamora. In \cite{Cam12}, Camere showed that for a complex K3 surface $X$ and an ample and globally generated line bundle $L$ on $X$, $M_L$ is stable with respect to $L$. Moreover, she showed that the same conclusion holds when $X$ is a complex abelian surface, $L$ is an ample and globally generated line bundle with $L^2\geq 14$ (\cite{Cam12}, Theorem 1 and Theorem 2). In \cite{ELM13}, Ein, Lazarsfeld and Mustopa showed that given a polarized smooth projective surface $(X,L)$ over an algebraically closed field of arbitrary characteristic and an arbitrary divisor $P$, $M_{dL+P}$ is stable with respect to $L$ for $d>>0$ (\cite{ELM13}, Theorem A). Further, they showed that $M_{dL}$ is cohomologically stable with respect to $L$ for polarized smooth projective varieties $(X,L)$ where $\dim(X)\geq 2$ and $d>>0$ provided $\textrm{Pic}(X)=\mathbb{Z}[A]$ (\cite{ELM13}, Proposition C). 

Recently, Torres-L\'opez and Zamora have showed that for a complex regular projective surface $X$, and an ample and globally generated line bundle $L$ on $X$, $M_L$ is stable with respect to $L$ if $-L\cdot K_X\geq 2$ (\cite{TLZ21}, Proposition 3.2). It follows that for anticanonical complex rational surfaces $X$, $M_L$ is stable with respect to $L$ if $L$ is ample and globally generated (see also \cite{TLZ21}, Corollaries 3.3, 3.4) since $-L\cdot K_X\geq 2$ for these surfaces (see for example (\cite{Har97}, Theorem III.1). Moreover, they showed that the same conclusion holds when $X$ is a complex Enriques surface and there exists a smooth irreducible curve $C\in|L|$ such that Cliff$(C)\geq 2$ (\cite{TLZ21}, Corollary 3.5).

It follows from a recent result of Caucci and Lahoz that $M_L$ is cohomologically stable with respect to $L$ for any ample and globally generated line bundle $L$ on any abelian surface over an algebraically closed field of arbitrary characteristic (\cite{CL21}, Theorem 1.5). Further, they proved effective semistability results for abelian varieties of arbitrary dimensions. Our main result of this note is the following.

\begin{theorem}
($=$ Theorem \ref{mainthm}) Let $X$ be an Enriques (resp. bielliptic) surface over an algebraically closed field of characteristic $\neq 2$ (resp. $\neq 2,3$). Then $M_L$ is cohomologically stable with respect to $L$ for any ample and globally generated line bundle $L$ on $X$.
\end{theorem} 

In particular, Theorem \ref{mainthm} improves \cite{TLZ21}, Corollary 3.5. Since cohomological stability implies stability, combining our result with the works of Camere, and Caucci--Lahoz, we obtain the following

\begin{theorem}
($=$ Corollary \ref{mainthm'}) Let $X$ be a smooth minimal complex projective surface of Kodaira dimension zero. Then $M_L$ is stable with respect to $L$ for any ample and globally generated line bundle $L$ on $X$.
\end{theorem}  

In the proof of our results, we combine the argument of Coand\u{a} (see \cite{Coa11}, \cite{ELM13}, proof of Proposition C, and \cite{CL21}, proof of Theorem 1.5) with the structure of the movable line bundles on Enriques and bielliptic surfaces.

By \cite{Cam12}, Proposition 2, for smooth regular complex projective surface with an ample and globally generated line bundle $L$, if the multiplication map 
\begin{equation}\label{mult}
    H^0(K_X)\otimes H^0(L)\to H^0(K_X\otimes L)
\end{equation}
surjects then $M_L$ is rigid. Consequently, for such an $L$, $M_L$ corresponds to an isolated point in the moduli of stable bundles with respect to $L$ if $X$ is a complex K3 or abelian surface. However, note that the map in \eqref{mult} can not surject if $X$ is a complex Enriques or bielliptic surface. In Proposition \ref{unobs}, we show that in the case of Enriques surfaces when char$(\mathbb{k})\neq 2$ and in the case of bielliptic surfaces when char$(\mathbb{k})\neq 2,3$, the moduli is smooth at $[M_L]$ and we compute the the dimension of its moduli component.

\vspace{5pt}

\noindent\textbf{Conventions.} We work over a fixed algebraically closed field $\mathbb{k}$. For a smooth projective variety $X$, $K_X$ denotes its canonical bundle. We shall denote $\chi(\mathcal{O}_X)$ by $\chi(X)$. We use the additive and multiplicative notation interchangeably for tensor product of line bundles, and $L^*$ denotes the dual of the line bundle $L$. Further, we shall use $L^{r}$ to denote the intersection product. The sign ``$\equiv$'' will be used for numerical equivalence and the sign ``$\sim$'' will be used for linear equivalence. 

\vspace{5pt}

\noindent\textbf{Acknowledgements.} We are grateful to Professor Purnaprajna Bangere and Professor Yusuf Mustopa for very helpful comments and suggestions on the article. We thank Professor A. G. Zamora for sharing the final version of \cite{TLZ21} with us. The first author was supported by the National Science Foundation, Grant No. DMS-1929284 while in residence at the Institute for Computational and Experimental Research in Mathematics in Providence, RI, as part of the ICERM Bridge program. The research of the second author was supported by a Simons Postdoctoral Fellowship from the Fields Institute for Research in Mathematical Sciences. We are also grateful to the anonymous referee for the suggestions and corrections that helped in the improvement of the exposition.

\section{Preliminaries on Enriques and bielliptic surfaces} 

 In this section we recall a few facts about Enriques and bielliptic surfaces that we shall require in the proof of our main theorem.

\subsection{Enriques surfaces} We start with the definition of Enriques surfaces.

\begin{definition}
(\cite{Cos83}) An {\it Enriques surface} $X$ over an algebraically closed field of characteristic $\neq 2$ is a smooth minimal projective surface with $q(X):=h^1(\mathcal{O}_X)=0$, and $K_X\neq 0$ but $2K_X=0$ (thus $p_g(X):=h^2(\mathcal{O}_X)=0$ and $\chi(X)=1$).
\end{definition}

We also need the following Bertini type theorem of Cossec on Enriques surfaces, which is valid over any algebraically closed field of characteristic $\neq 2$.

\begin{theorem}\label{cossec}
(\cite{Cos83} Theorem 1.5.1 and Remark 1.6.3. (i)) Assume char$(\mathbb{k})\neq 2$. Let $X$ be an Enriques surface and let $L$ be a movable (i.e. $|L|$ is non-empty and fixed component free) line bundle on $X$.
\begin{itemize}
    \item[(1)] If $L^2>0$ then $h^0(L)=\chi(L)=1+\frac{1}{2}L^2$.
    \item[(2)] If $L^2=0$, then $L\sim kP$ where $|P|$ is an elliptic pencil on $X$ and $h^0(L)=k+1$ ($k\geq 1$). Moreover, every elliptic pencil $|P|$ on $X$ has exactly two multiple fibres. These are double fibres $2E_1$ and $2E_2$ where the divisors $E_i$ is a smooth elliptic curve, a rational curve with an ordinary double point, or a loop of $b\geq 2$ nodal curves. 
\end{itemize}
\end{theorem}

\subsection{Bielliptic surfaces} We now move on to bielliptic surfaces and define them. Bielliptic surfaces in positive characteristic have recently been studied by Boada De Narvaez (\cite{BDN}). 

\begin{definition}
(\cite{BDN}) A {\it bielliptic surface} $X$ over an algebraically closed field of characteristic $\neq 2,3$ is a smooth minimal projective surface with Kodaira dimension zero that satisfies $q(X):=h^1(\mathcal{O}_X)=1$ and $p_g(X):=h^2(\mathcal{O}_X)=0$ (thus $\chi(X)=0$).
\end{definition} 

Recall that a surjective morphism $f:Y\to B$ from a smooth projective surface to a smooth projective curve is called a {\it genus one fibration} if $f_*\mathcal{O}_Y=\mathcal{O}_B$ and all but finitely many fibres of $f$ are integral curves of arithmetic genus one. Further, it is called {\it minimal} if no fibre of $f$ contains any $(-1)$ curve. Also, a minimal genus one fibration is called {\it elliptic} if the general fibres are smooth. We now recall a characterization of these surfaces. 
\begin{proposition}
(\cite{Bea}, Chapter VI, \cite{BDN}, Chapter 1) Assume char$(\mathbb{k})\neq 2,3$. Given a bielliptic surface $X$, there exists two elliptic curves $A$ and $B$, and an abelian group $G$ acting on $A$ and on $B$ such that 
\begin{itemize}
    \item[(1)] $A/G$ is elliptic, and $B/G\cong\mathbb{P}^1$, 
    \item[(2)] $S\cong (A\times B)/G$ where $G$ acts on $A\times B$ componentwise, and
    \item[(3)] there are two elliptic fibrations $\varphi:X\to B/G\cong\mathbb{P}^1$ and $\psi:X\to A/G$.
\end{itemize}
\end{proposition}
Throughout this short note, we shall make use of the following diagram that exists by the proposition above.
\[
\begin{tikzcd}
& X\cong (A\times B)/G \arrow[dl, swap, "\varphi"] \arrow[dr, "\psi"] & \\
B/G\cong \mathbb{P}^1  & & A/G
\end{tikzcd}
\]
The morphism $\psi$ is smooth with every fibre isomorphic to $B$. All smooth fibres of $\varphi$ are isomorphic to $A$ and it has finitely many multiple fibres with multiplicities $m_1,\cdots, m_t$ ($t\geq 1$ is the number of multiple fibres) all of which are multiples of smooth elliptic curves. We shall denote by $A$ and $B$, the class in Num$(X)$ of the fibres of $\varphi$ and $\psi$ respectively. The classes intersect as $$A^2=0,\quad B^2=0,\quad A\cdot B=\gamma:=|G|.$$ 
If one denotes the $\textrm{lcm}(m_1,\cdots,m_t)$ by $\mu$, then it follows from \cite{Ser}, Theorem 1.4 in characteristic zero, and by \cite{BDN}, Theorem 3.3) in characteristic $\neq 2,3$ that a basis of $\textrm{Num}(X)$ is given by $\left\{A/\mu,(\mu/\gamma) B\right\}$. Also, by \cite{Ser}, Lemma 1.3, and \cite{BDN}, Lemma 3.2, a line bundle $L\equiv a\cdot A/\mu+b\cdot (\mu/\gamma)B$ is effective only if $a,b\geq 0$, ample if and only if $a,b>0$.\par 

\smallskip

We aim to prove an analogue of Theorem ~\ref{cossec} for bielliptic surfaces. The proof of the Lemma below is given by Mella and Palleschi in \cite{MP93}, Proposition 2.4, p. 223 in characteristic zero. The exact same proof goes through for arbitrary characteristics and we reproduce it here.

\begin{lemma}\label{fixedlemma}
Assume char$(\mathbb{k})\neq 2,3$. If $D=\sum_{i}b_iA_i$ with $0\leq b_i<m_i$ is an effective divisor where $A_i$ is the reduced part of the $i$-th multiple fibre, then $h^0(\mathcal{O}_X(D))=1$.
\end{lemma}

\noindent\textit{Proof.} Consider $S:=\left\{D=\sum_{i}b_iA_i|0\leq b_i<m_i\textrm{ and $b_j\neq0$ for some $j$}\right\}$. Since $S$ is finite, if there exists $D\in S$ with $h^0(\mathcal{O}_X(D))\geq 2$, then there exists $D'\in |\mathcal{O}_X(D)|$ such that $D'=A+\sum_{i}e_iA_i$ with $e_i\geq 0$ since support of any $D'\in |\mathcal{O}_X(D)|$ is contained in a finite union of fibres of $\varphi$. Set $$b_0:=\textrm{min}\left\{b\mid\textrm{ there exists } D\in S \textrm{ with }D\equiv bA/\mu\textrm{ and }h^0(D)\geq 2\right\}$$
and $D_0\in S$ is the divisor for which the minimum is attained. In what follows $b_0\geq \mu$. Let $D_0=\sum_ic_iA_i$ with $0\leq c_i<m_i$ with $c_1\neq 0$. By the discussion above, there exists $D_0'=A+\sum_id_iA_i$ with $d_i\geq 0$. Observe that $D_0\neq c_1A_1$ since $b_0\geq \mu$, and thus $D_0-c_1A_1\in S$ and $h^0(\mathcal{O}_X(D_0-c_1A_1)\geq 2$ that contradicts the minimality of $b_0$.\QEDB\par

\vspace{5pt}

Now we prove an analogue of Theorem ~\ref{cossec}, the main argument of the proof of Proposition \ref{bertinibi} (b) is again taken from \cite{MP93}, Proposition 2.4 (the only observation being that the authors did not require the hypothesis of $0$--very ampleness, instead the weaker hypothesis of the line bundle being movable suffices).

\begin{proposition}\label{bertinibi}
Let $X$ be a bielliptic surface and $L$ be a movable line bundle on $X$. Assume char$(\mathbb{k})\neq 2,3$.
\begin{itemize}
    \item[(a)] If $L^2>0$ then $h^0(L)=\chi(L)=\frac{1}{2}L^2$. 
    \item[(b)] If $L^2=0$, then exactly one of the following cases happens.
\begin{itemize}
    \item[(1)] $L\cdot A=0$, $L\sim \varphi^*\mathcal{O}_{\mathbb{P}^1}(a)$ and $h^0(L)=a+1$ for some $a\geq 1$.
    \item[(2)] $L\cdot B=0$, $L\sim \psi^*L'$ where $L'$ is a line bundle on $A/G$ with $\textrm{deg}(L')\geq 2$. In this case, we have $h^0(L)=\textrm{deg}(L')$.
\end{itemize}
\end{itemize}
\end{proposition}

\noindent\textit{Proof.} (a) Since $\textrm{Num}(X)=\mathbb{Z}[A/\mu]\oplus \mathbb{Z}[(\mu/\gamma)B]$, it follows that if $L\equiv aA/\mu+b(\mu/\gamma)B$ with $a,b\in \mathbb{Z}$, then $a,b>0$. Thus we obtain that $L$ is ample. The assertion follows from \cite{Ser}, Lemma 1.3 (iii) for characteristic zero, and from \cite{BDN}, Lemma 3.2 (iv) for characteristic $\neq 2,3$.

(b) Since $\textrm{Num}(X)=\mathbb{Z}[A/\mu]\oplus \mathbb{Z}[(\mu/\gamma)B]$, hence if $L^2=0$, then either $L\cdot A=0$ or $L\cdot B=0$. \par 
\smallskip
\noindent\underline{Case 1: $L\cdot A=0$.} In this case, the support of an effective divisor is contained in a finite union of fibres of $\varphi$. Notice that two smooth fibres of $\varphi$ are linearly equivalent. For the sake of contradiction, suppose there exists $D_0\in |L|$ that can be written as follows where $A_i$ is the reduced part of the $i$-th multiple fibre:
\begin{equation}\label{eq1}
    D_0= aA+\sum_{i=1}^{n}a_iA_i,\qquad a\geq 0,\, 0\leq a_i\leq m_i-1,\, a_1\neq 0.
\end{equation}
Also, since $L$ is movable, there exists $D_1\in |L|$ that can be expressed as follows:
\begin{equation}\label{eq2}
    D_1= bA+\sum_{i=2}^{n}b_iA_i,\qquad b\geq 0,\, 0\leq b_i\leq m_i-1.
\end{equation}
Now, ~\eqref{eq1} and ~\eqref{eq2} shows that we have the following
\begin{equation}\label{eq3}
    (a-b)A+\sum_{i=1}^n c_iA_i\sim \sum_{i=2}^n d_iA_i,\qquad c_1=a_1,\, c_id_i=0\,\textrm{ and }\, 0\leq c_i,d_i\leq m_i-1\textrm{ for all}\, i\neq 2.
\end{equation}
By setting $D= \sum_{i=1}^n c_iA_i$ if $a-b<0$, and $\sum_{i=2}^n d_iA_i$ otherwise, we see from ~\eqref{eq3} that $h^0(\mathcal{O}_X(D))\geq 2$ that contradicts Lemma \ref{fixedlemma}. 
Consequently, since $\varphi_*\mathcal{O}_X=\mathcal{O}_{\mathbb{P}^1}$, 
we obtain $h^0(L)=h^0(\mathcal{O}_{\mathbb{P}^1}(a))=a+1$ (clearly $a\neq 0$ since $h^0(L)\geq 2$).\par 
\smallskip
\noindent\underline{Case 2: $L\cdot B=0$.} In this case, the fibres of an effective divisor of $|L|$ is contained inside the fibres of $\psi$. Since the fibres are smooth, we get $L\sim\psi^*L'$ where $\textrm{deg}(L')\geq 0$ and as before $h^0(L)=h^0(L')=\textrm{deg}(L')$. Clearly $\textrm{deg}(L')\geq 2$ since $h^0(L)\geq 2$.\QEDB\par 

\section{Syzygy bundles on Enriques and bielliptic surfaces}

In this section, we show that the syzygy bundles $M_L$ are stable on Enriques (resp. bielliptic) surfaces $X$ when char$(\mathbb{k})\neq 2$ (resp. char$(\mathbb{k})\neq 2,3$) for ample and globally generated line bundles $L$. Further, we show that over complex such surfaces, they are smooth points of the moduli of stable vector bundles and we compute the dimension of the components containing them. Observe that it follows from \eqref{*} that $c_1(M_L)=L^*$ and $c_2(M_L)=L^2$.

\subsection{Stability of syzygy bundles} We start with the definition and a result on slope stability. We then recall a vanishing theorem of Green that was used in \cite{Coa11}, \cite{ELM13} and in \cite{CL21}.

Let $(X,L)$ be a smooth polarized projective variety of dimension $n$, and let $E$ be a torsion-free coherent sheaf on $X$. The {\it slope} of $E$ with respect to $L$ is by definition $$\mu_L(E):=\frac{c_1(E)\cdot L^{n-1}}{\textrm{rank}(E)}.$$ 

\begin{definition}
Let $(X,L)$ be a smooth polarized projective variety of dimension $n$, and let $E$ be a vector bundle on $X$. $E$ is called {\it stable with respect to $L$} if for all subsheaf $F\subsetneq E$ with $0<\textrm{rank}(F)<\textrm{rank}(E)$, the inequality $\mu_L(F)<\mu_L(E)$ holds.
\end{definition}

\begin{lemma}\label{strongs}
Let $(X,L)$ be a smooth polarized projective variety and let $E$ be a vector bundle on $X$. If for any $0<r<\textrm{rank}(E)$ and for any line bundle $N$ on $X$ with $\mu_L(\wedge^{r}E\otimes N)\leq 0$, one has $$H^0(\wedge^{r}E\otimes N)=0,$$ then $E$ is stable with respect to $L$.
\end{lemma}

If $E$ satisfies the condition of the above theorem, then it is called {\it cohomologically stable} with respect to $L$. The following vanishing result is due to Green.

\begin{lemma}\label{green}
(\cite{Gre84}, Theorem (3.a.1), \cite{Coa11}, Lemma 1.2) Let $L$ be a globally generated line bundle on a smooth projective variety $X$. Then, for a line bundle $N$ on $X$, $H^0(\wedge^r M_L\otimes N)=0$ if $r\geq h^0(N)$.
\end{lemma}

We now state and prove a lemma that we shall use in the proof of our main theorem. 
\begin{lemma}\label{mainlemma}
Let $X$ be an Enriques (resp. bielliptic) surface and we assume char$(\mathbb{k})\neq 2$ (resp. char$(\mathbb{k})\neq 2,3$). Let $L$ be an ample and globally generated line bundle on $X$. Moreover, let $M$ be a movable line bundle on $X$ with $M^2>0$. If $\mu_L(\wedge^r M_L\otimes M)\leq 0$ for some integer $0<r<\textrm{rank}(M_L)=h^0(L)-1$, then $r\geq h^0(M)$.
\end{lemma}
\noindent\textit{Proof.} Notice that $\mu_L(\wedge^rM_L\otimes M)=-\frac{rL^2}{h^0(L)-1}+L\cdot M$. Thus, $\mu_L(\wedge^r M_L\otimes M)\leq 0$ is equivalent to the following
\begin{equation}\label{maineq}
    r\geq (L\cdot M)\frac{h^0(L)-1}{L^2}.
\end{equation}
By Hodge index theorem, $(L\cdot M)^2\geq L^2\cdot M^2$. Suppose $M^2\geq L^2$, then by ~\eqref{maineq} 
\begin{equation*}
    r\geq \sqrt{L^2}\sqrt{M^2}\frac{h^0(L)-1}{L^2}\geq h^0(L)-1 
\end{equation*}
that leads to a contradiction. Thus, we conclude that $L^2>M^2$. Further, notice that by Theorem \ref{cossec} and Proposition \ref{bertinibi}, we have $h^0(M)=\chi(M)=\chi(X)+\frac{1}{2}M^2$ and $h^0(L)=\chi(L)=\chi(X)+\frac{1}{2}L^2$. Consequently, Hodge index theorem and ~\eqref{maineq} yields
\begin{equation}\label{maineq2}
    r> M^2\frac{h^0(L)-1}{L^2}\geq 2(h^0(M)-\chi(X))\left(\frac{1}{2}+\frac{\chi(X)-1}{L^2}\right).
\end{equation}
Thus, when $\chi(X)=1$ (i.e. $X$ is an Enriques surface), we obtain $r>h^0(M)-1$ and the conclusion follows. On the other hand, when $\chi(X)=0$ (i.e. $X$ is a bielliptic surface), we obtain $r>h^0(M)-\frac{2h^0(M)}{L^2}$. Notice that $2h^0(M)=M^2$ and thus $\frac{2h^0(M)}{L^2}<1$ since we already know that $L^2>M^2$. Thus, we obtain $r>h^0(M)-1$ and the proof is now complete.\QEDB\par

\vspace{5pt}

Now we are ready to prove our main result.

\begin{theorem}\label{mainthm}
Let $X$ be an Enriques (resp. bielliptic) surface and we assume char$(\mathbb{k})\neq 2$ (resp. char$(\mathbb{k})\neq 2,3$). Let $L$ be an ample and globally generated line bundle on $X$. Then the syzygy bundle $M_L$ is cohomologically stable with respect to $L$.
\end{theorem}

\noindent\textit{Proof.} We apply Lemma ~\ref{strongs}: we assume that $\mu_L(\wedge^r M_L\otimes N)\leq 0$ for a line bundle $N$ and an integer $0<r<\textrm{rank}(M_L)=h^0(L)-1$ and we aim to show that $H^0(\wedge^r M_L\otimes N)=0$. By Lemma ~\ref{green}, it is enough to show that $r\geq h^0(N)$. Thus, we may assume that $h^0(N)\geq 2$. As in the proof of Lemma ~\ref{mainlemma}, the inequality $\mu_L(\wedge^r M_L\otimes N)\leq 0$ is equivalent to 
\begin{equation}\label{befored}
    r\geq (L\cdot N)\frac{h^0(L)-1}{L^2}.
\end{equation}
Now we decompose $N$ as $N=M+F$ where $M$ is movable and $F$ is fixed. Notice that $h^0(M)=h^0(N)$. Since $L$ is ample, it follows from ~\eqref{befored} that 
\begin{equation}\label{afterd}
    r\geq (L\cdot M)\frac{h^0(L)-1}{L^2},
\end{equation}
and we aim to show that $r\geq h^0(M)$. Notice that by ~\eqref{afterd} $L\cdot M<L^2$ since $0<r<h^0(L)-1$. The required inequality immediately follows from Lemma ~\ref{mainlemma} when $M^2>0$. Thus, we assume that $M^2=0$.
\smallskip

\noindent\underline{Case 1: $X$ is an Enriques surface.} In this case, the inequality ~\eqref{afterd} boils down to $r\geq \frac{1}{2}(L\cdot M)$ since by Theorem \ref{cossec} $h^0(L)=\chi(L)$. Also, by Theorem ~\ref{cossec}, we know that $|M|=|kP|$ where $|P|$ is an elliptic pencil and $h^0(M)=k+1$ where $k\geq 1$. We claim that $L\cdot M\geq 4k$. To see this, we use the structure of the elliptic pencils on $X$ (i.e. Theorem ~\ref{cossec}). Notice that $L\cdot M=kL\cdot P$. Since $|P|$ contains a multiple fibre $2E$, one obtains $L\cdot M=2kL\cdot E$. Since $L$ is ample, the claim follows if $E$ is a loop of $b\geq 2$ nodal curves. Thus, we may assume that $E$ is irreducible, i.e., a smooth elliptic curve or a rational curve with an ordinary double point. In both cases, since the arithmetic genus $p_a(E)=1$, it follows that $L\cdot E\geq 2$ since $L|_E$ is globally generated as $L$ is globally generated and that proves the claim. Consequently, $r\geq 2k\geq k+1$ since $k\geq 1$. 
\smallskip

\noindent\underline{Case 2: $X$ is a bielliptic surface.} In this case, 
~\eqref{afterd} becomes 
\begin{equation}\label{again}
    r\geq \frac{1}{2}(L\cdot M)-\frac{L\cdot M}{L^2}.
\end{equation}
We are going to have two cases by Proposition ~\ref{bertinibi}. First, assume $M\cdot A=0$. We know that $M\sim \varphi^*\mathcal{O}_{\mathbb{P}^1}(a)$ with $a\geq 1$ and $h^0(M)=a+1$. Since $\varphi$ has multiple fibres, it follows that $L\cdot M\geq 2a L\cdot A_i$ for any $i\in\{1,\cdots,t\}$ where $A_i$'s are the reduced part of the multiple fibres. Since $A_i$ is a smooth elliptic curve and $L$ is ample and globally generated, as before we have $L\cdot A_i\geq 2$. Consequently, $L\cdot M\geq 4a$. We obtain from ~\eqref{again} that $$r\geq 2a-\frac{L\cdot M}{L^2}>2a-1\implies r\geq 2a.$$ The conclusion follows since $2a\geq a+1$ as $a\geq 1$. Next, assume $M\cdot B=0$. In this case, we know that $M\sim\psi^*M'$ and $h^0(M)=\textrm{deg}(M')\geq 2$. Since support of a divisor in the linear series of $|M|$ is a union of $\textrm{deg}(M')$ fibres of $\psi$, and $L\cdot E\geq 2$ for any smooth elliptic curve $E$ on $X$ (as $L$ is ample and globally generated), it follows that $L\cdot M\geq 2\textrm{deg}(M')$. Thus, by ~\eqref{again}
$$r\geq \textrm{deg}(M')-\frac{L\cdot M}{L^2}>\textrm{deg}(M')-1.$$
That completes the proof.\QEDB\par 

\begin{corollary}\label{mainthm'}
Let $X$ be a smooth minimal complex projective surface of Kodaira dimension zero. Then $M_L$ is stable with respect to $L$ for any ample and globally generated line bundle $L$ on $X$.
\end{corollary}

\noindent\textit{Proof.} Recall that a smooth minimal complex projective surface of Kodaira dimension zero is either a K3, or an abelian, or an Enriques, or a bielliptic surface. If $X$ is a K3 surface, then the conclusion follows from \cite{Cam12}, Theorem 1 (also by \cite{TLZ21}, Corollary 3.5). If $X$ is either an Enriques or a bielliptic surface, then the assertion follows from Theorem \ref{mainthm}. Finally, if $X$ is an abelian surface, then notice that $L^2\geq 2$ and $L\cdot E\geq 2$ for any elliptic curve $E$ on $X$ since $L$ is ample and globally generated. Consequently, the conclusion follows from \cite{CL21}, Theorem 1.5.\QEDB\par

\subsection{Moduli points of syzygy bundles} Let $(X,L)$ be a smooth polarized projective surface and let $E$ be a vector bundle with $\textrm{rank}(E)=r$ and $c_i:=c_i(E)$ for $i=1,2$, stable with respect to $L$. Let us denote the moduli space of $L$-stable vector bundles on $X$ with rank $r$ and Chern classes $c_i$ by $M_{X,L}:=M_{X,L}(r,c_1,c_2)$. Recall from \cite{HL} that if $\textrm{Ext}^2(E,E)=0$, then $M_{X,L}$ is smooth at $[E]$ and its dimension at $[E]$ is 
\begin{equation}\label{dim}
    2rc_2-(r-1)c_1^2-(r^2-1)\chi(X).
\end{equation} 
Since the syzygy bundle $M_L$ is $L$-stable for an ample and globally generated line bundle $L$ on Enriques or bielliptic surfaces $X$ by Theorem \ref{mainthm}, they are points of $M_{X,L}$. In the following proposition we show the smoothness of the moduli at the syzygy bundles and compute the dimension of the component containing them.

\begin{proposition}\label{unobs}
Let $X$ be a smooth projective Enriques (resp. bielliptic) surface and assume char$(\mathbb{k})\neq 2$ (resp. char$(\mathbb{k})\neq 2,3$). Let $L$ be an ample and globally generated line bundle on $X$. Then $\textrm{Ext}^2(M_L,M_L)=0$. In particular, $[M_L]$ is a smooth point of $M_{X,L}$ and the dimension of the component containing $[M_L]$ is $\frac{1}{4}(L^2)^2+L^2+1$ (resp. $\frac{1}{2}(L^2)^2$).
\end{proposition}

\noindent\textit{Proof.} Tensoring \eqref{*} by $K_X$ and taking global sections, one obtains $H^0(M_L\otimes K_X)=0$ thanks to $p_g(X)=0$. We claim that $H^0(M_L\otimes M_L^*\otimes K_X)=0$. To see this, dualize \eqref{*} to obtain 
\begin{equation}\label{**}
    0\to L^*\to H^0(L)^*\otimes \mathcal{O}_X\to M_L^*\to 0.
\end{equation}
Tensoring \eqref{**} by $M_L\otimes K_X$ and taking the long exact sequence of cohomology, it is easy to see that it is enough to prove that $H^1(M_L\otimes L^*\otimes K_X)=0$ to prove the claim. To this end, tensoring \eqref{*} by $K_X\otimes L^*$ and taking the long exact sequence of cohomology the required vanishing follows thanks to $p_g(X)=0$ and $H^1(L)=0$ (by Theorem \ref{cossec}, and Proposition \ref{bertinibi}). Thus $\dim\textrm{Ext}^2(M_L,M_L)=h^2(M_L\otimes M_L^*)=h^0(M_L\otimes M_L^*\otimes K_X)=0$ where the second equality is obtained by duality. Finally, the dimension of the moduli component follows from \eqref{dim}.\QEDB\par

\bibliographystyle{plain}

\begin{thebibliography}{100}

\bibitem[Boa21]{BDN}
  Boada De Narvaez, Daniel.
  \emph{Moduli of Bielliptic Surfaces}.
  Thesis (Ph.D.) Technische Universit\"at M\"unchen, 2021. 

\bibitem[Bea83]{Bea}
  Beauville, Arnaud.
  \emph{Complex algebraic surfaces}.
  Translated from the French by R. Barlow, N. I. Shepherd-Barron and M. Reid. London Mathematical Society Lecture Note Series, 68. Cambridge University Press, Cambridge, 1983. {\rm iv}+132 pp.

\bibitem[Bre08]{Bre08}  
   Brenner, Holger
  \emph{Looking out for stable syzygy bundles}.
   With an appendix by Georg Hein. Adv. Math. 219 (2008), no. 2, 401-427.  

\bibitem[Cam12]{Cam12}  
   Camere, Chiara.
  \emph{About the stability of the tangent bundle of $\mathbb{P}^n$ restricted to a surface}.
   Math. Z. 271 (2012), no. 1-2, 499--507.
   
\bibitem[Cos83]{Cos83}  
   Cossec, Fran\c{c}ois R.
  \emph{Projective models of Enriques surfaces}.
   Math. Ann. 265 (1983), no. 3, 283–334.   
   
   
\bibitem[CL21]{CL21}  
   Caucci, Federico; Lahoz, Mart\'i
  \emph{Stability of syzygy bundles on abelian varieties}.
   Preprint, available at \url{https://arxiv.org/abs/2007.08846}.  
   
\bibitem[CMM10]{CMM10}  
   Costa, Laura; Macias Marques, Pedro; Mir\'o-Roig, Rosa Mar\'ia.
  \emph{Stability and unobstructedness of syzygy bundles}.
    J. Pure Appl. Algebra 214 (2010), no. 7, 1241–1262

\bibitem[Coa11]{Coa11}  
   Coand\u{a}, Iustin.
  \emph{On the stability of syzygy bundles}.
   Internat. J. Math. 22 (2011), no. 4, 515--534.
   
\bibitem[EL92]{EL92}  
   Ein, Lawrence; Lazarsfeld, Robert.
  \emph{Stability and restrictions of Picard bundles, with an application to the normal bundles of elliptic curves}.
   Complex projective geometry (Trieste, 1989/Bergen, 1989), 149–156, London Math. Soc. Lecture Note Ser., 179, Cambridge Univ. Press, Cambridge, 1992.   
   
\bibitem[ELM13]{ELM13}  
   Ein, Lawrence; Lazarsfeld, Robert; Mustopa, Yusuf.
  \emph{Stability of syzygy bundles on an algebraic surface}.
   Math. Res. Lett. 20 (2013), no. 1, 73--80.

\bibitem[Fle84]{Fle84}  
   Flenner, Hubert.
  \emph{Restrictions of semistable bundles on projective varieties}.
   Comment. Math. Helv. 59 (1984), no. 4, 635–650.

\bibitem[Gre84]{Gre84}  
   Green, Mark L.
  \emph{Koszul cohomology and the geometry of projective varieties}.
   J. Differential Geom. 19 (1984), no. 1, 125--171.
   
\bibitem[Har97]{Har97}  
   Harbourne, Brian.
  \emph{Anticanonical rational surfaces}.
   Trans. Amer. Math. Soc. 349 (1997), no. 3, 1191–1208.  
   
\bibitem[HL97]{HL}  
   Huybrechts, Daniel; Lehn, Manfred.
  \emph{The geometry of moduli spaces of sheaves}.
   Aspects of Mathematics, E31. Friedr. Vieweg \& Sohn, Braunschweig, 1997. {\rm xiv}+269 pp.  
   
\bibitem[MM11]{MM11}  
   Macias Marques, Pedro; Mir\'o-Roig, Rosa Mar\'ia.
  \emph{Stability of syzygy bundles}.
    Proc. Amer. Math. Soc. 139 (2011), no. 9, 3155--3170.    

\bibitem[MP93]{MP93}  
   Mella, M.; Palleschi, M.
  \emph{The k-very ampleness on an elliptic quasi bundle}.
    Abh. Math. Sem. Univ. Hamburg 63 (1993), 215–226.
    
\bibitem[Ser90]{Ser}
  Serrano, Fernando.
  \emph{Divisors of bielliptic surfaces and embeddings in ${\bf P}^4$}.
  Math. Z. 203 (1990), no. 3, 527--533.     
   
\bibitem[TLZ21]{TLZ21}  
  Torres-L\'{o}pez, H.; Zamora, A. G.
  \emph{Some remarks on H-stability of syzygy bundle on algebraic surface}.
  Beitr\"age zur Algebra und Geometrie (2021).  
\bibitem[Tri10]{Tri10}  
  Trivedi, V.
  \emph{Semistability of syzygy bundles on projective spaces in positive characteristics}.
  Internat. J. Math. 21 (2010), no. 11, 1475--1504. 
   
 
  
   






  

\end{thebibliography}

\end{document}